\documentclass[12pt, a4paper]{article}
\usepackage[margin=25.4mm]{geometry}
\usepackage[T1]{fontenc}
\usepackage{newtxtext, newtxmath}
\usepackage[ruled]{algorithm2e}
\usepackage{titlesec}
\titleformat*{\section}{\raggedright\large\bfseries}
\titleformat*{\subsection}{\raggedright\normalsize\bfseries}
\usepackage[colorlinks]{hyperref}
\usepackage[numbers, sort&compress]{natbib}
\usepackage{paralist}
\usepackage[normalem]{ulem}
\renewcommand*{\underline}[2]{\uline{#1} (Origin: #2)}
\usepackage[english]{babel}                     
\usepackage{blindtext}
\usepackage{setspace}
\usepackage{lscape}
\usepackage{array}
\usepackage{hyperref}
\usepackage{appendix}
\usepackage{autobreak}
\hypersetup{hidelinks} 
\usepackage{ breqn}
\usepackage{autobreak}
\usepackage{graphicx}                           
\usepackage{amsmath}                            
\usepackage{hyperref}
\usepackage{algorithmicx}
\usepackage{lastpage}                           
\usepackage{lipsum}                             
\usepackage{indentfirst}
\usepackage{amsfonts} 
\usepackage{amsopn}
\usepackage{array}
\allowdisplaybreaks
\begin{document}
	\begin{flushleft}\Large\textbf{
			Four limit cycles in three-dimensional Lotka-Volterra competitive systems with classes 28, 30 and 31 in Zeeman’s classification via automatic search}
		
	\end{flushleft}
	\begin{flushleft}
		\textbf{Mingzhi Hu$^{1}$, Zhengyi Lu$^{1}$ and Yong Luo$^{2}$}\\[1em]
		
			$^{1}$ School of Mathematical Sciences, Sichuan Normal University,  Chengdu 610068, China
			
			$^{2}$ College of Mathematics and Physics, Wenzhou University,  Wenzhou 325000, China
		
		\hspace{1em}
		
		Correspondence should be addressed to Yong Luo; luoyong@wzu.edu.cn
	\end{flushleft}
	
	\section*{Abstract}
	Four limit cycles are constructed for classes 28, 30 and 31 in Zeeman's classification, together with the results in \cite{GY2009} for class 27, \cite{w} for class 29 and \cite{YP} for class 26 which indicate that for each class among classes 26$-$31, there exist at least four limit cycles. This gives a partial answer to a problem proposed in \cite{Ho} as well as in \cite{YP}.
	
	\vspace{1em}
	
	\noindent Keywords: Lotka-Volterra, competitive systems, Zeeman's classification, four limit cycles
	\vspace{1em}
	
	\section{Introduction}
	
	In addition to the Liapunov function originally given by Volterra \cite{Vol}, MacArthur \cite{Mc} gave a Liapunov function for a  Lotka-Volterra system with a symmetric interaction matrix. Using MacArthur's Liapunov function, it can be shown that a $\omega -$limit set of a competitive system with a symmetric interaction matrix is either an isolated equilibrium or a continuum of equilibria.
	
	In 1975, May and Leonard \cite{May} gave a landmark result for three-dimensional competitive systems. They constructed a three-dimensional Lotka-Volterra competitive system with a cyclic symmetric interaction matrix and obtained the periodicity of the system under certain conditions. Meanwhile, through numerical simulation, they found that a  three-dimensional cyclically symmetric competitive system may have a limit cycle and an attractive or a repulsive heteroclinic cycle.
	
	In 1979,  Coste, Peyraud and Coullet \cite{CPC} and Schuster, Sigmund and Wolff \cite{SSR} proved the existence of a limit cycle for three-dimensional competitive systems(which  actually belongs to class 27 in Zeeman's classification \cite{Zee}) independently and proved the attractiveness or repulsion of the heteroclinic cycle.  In 1981, Hofbauer \cite{Ho81} proved the existence of limit cycles of general $n -$ dimensional $(n \geqslant 3)$ competitive systems by Hopf bifurcation theory.
	
	Smale \cite{Smale} in  1976  and  Hirsch \cite{Hir88} in 1988 showed that the dynamical behavior of $n-$ dimensional competitive systems is comparable to that of $(n-1)-$ dimensional systems. For a three-dimensional Lotka-Volterra competitive system, Hirsch's result \cite{Hir88} ensured that there exists an invariant manifold (called a carrying simplex) that is homeomorphic to a two-dimensional simplex and attracts all orbits except the origin. Based on Hirsch's theorem, Zeeman \cite{Zee} used the qualitative analysis of the system at the boundary surface to define the combinatorial equivalence relation through parameter inequalities, and obtained thirty-three stable classes of  three-dimensional Lotka-Volterra competitive systems.  Zeeman's result indicated that the limit sets for 27 classes among the 33 ones are all fixed points \cite{Zee}, so the dynamical behaviors of systems in these 27 classes are fully described. Using Hopf bifurcation theorem, Zeeman further proved that among each remaining six classes: from class 26 to class 31, there is always a system can be constructed to have an isolated periodic orbit or limit cycles by selecting parameters.
	
	The first result of a three-dimensional Lotka-Volterra competitive system with two limit cycles was obtained in class 27 of Zeeman's classification by Hofbauer and So \cite{Ho} based on Hirsch's monotonic flow theorem, the center manifold theorem and the Hopf bifurcation theorem. In their example, a local stable positive equilibrium was surrounded by two limit cycles, one is from the Hopf bifurcation theorem and the other is guaranteed by the Poincaré-Bendixson theorem. In \cite{Ho}, Hofbauer and So proposed a problem of how many limit cycles may exist in the six classes 26$-$31 in Zeeman's classification.
	
	In 2002, Lu and Luo \cite{Lu02} constructed two limit cycles in class 26, class 28, and class 29, respectively. Hofbauer and So's questions about whether there are two limit cycles for other classes were partially answered.
	
	In 2003, the first attempt to construct three limit cycles was given by Lu and Luo \cite{Lu03}. They constructed three limit cycles in Zeeman's system of 27 classes. As Yu, Han and Xiao \cite{YP} pointed out that due to a sign error, the positive definiteness of the Liapunov function for calculating the focal value was not be guaranteed, therefore, Lu and Luo's attempt failed.
	
	In 2006, Gyllenberg, Yan and Wang \cite{GY2006} found that the system with two limit cycles in class 29 constructed by Lu and Luo \cite{Lu02} actually can have a third limit cycle. The instability of the outer small amplitude limit cycle together with the repulsion of the boundary of the carrying simplex ensure that the Poincaré-Bendixson theorem holds true which can ensure the existence of a third limit cycle.
	
	In 2008, Lian, Lu and Luo \cite{lian} found through automatic search that both class 30 and class 31 in Zeeman's classification can have systems with three limit cycles. Two of them are small-amplitude limit cycles, and the other is a large-scale limit cycle obtained by the Poincaré-Bendxison Theorem. 

	There were also works by Gyllenberg and Yan \cite{GY2007} in 2009 and Wang, Huang and Wu \cite{w} in 2011, respectively, for the existence of multiple limit cycles for class 30. Gyllenberg and Yan claimed to have constructed two limit cycles in class 30. One is a small-amplitude limit cycle obtained by the Hopf bifurcation, and the other is obtained by the Poincaré-Bendxison Theorem. As pointed out by Yu, Han and Xiao \cite{YP}, their construction of limit cycles was incorrect since the Liapunov function they gave is not positive definite. However, with careful calculation, three limit cycles can actually be constructed by using their original  example \cite{GY2007}. First, two small-amplitude limit cycles can be constructed such that the outer one is stable. Since the system belongs to class 30, the boundary of the simplex is   an attractor, therefore the third limit cycle can be obtained by using the Poincaré-Bendxison Theorem.  Three years later, Wang, Huang, and Wu \cite{w} proved that for each class in Zeeman's classification of classes 26$-$31 there exists a system with at least three small-amplitude limit cycles. Unfortunately, the result in their case of class 30 was incorrect, since the constructed system was not a competitive one.
	
	In 2009, Gyllenberg and Yan \cite{GY2009} constructed an example in class 27 similar to \cite{Lu03} and claimed that they obtained four limit cycles, three of which were obtained due to the Hopf bifurcation, and the fourth one was obtained through the existence of the heteroclinic cycle and using the Poincaré-Bendxison Theorem. Similarly, they ignored the positive definiteness of the Liapunov function. Yu, Han and Xiao \cite{YP} pointed out that in  \cite{GY2009}  there is an example belonging to class 27 which can have four limit cycles. A system in class 29 with four limit cycles was given by Wang, Huang, and Wu \cite{w}.
	
	Recently, Yu, Han, and Xiao \cite{YP} gave two examples  in class 27 and two in class 26 with four small-amplitude  limit cycles, respectively.
	
	In summary, multiple limit cycles can appear in each  six classes of  Zeeman's classifications : classes 26, 27, 28, 29, 30 and 31. The known results until now are as follows: for the classes 26, 27 and 29, there are at least 4 limit cycles.
	
	An open problem proposed by Yu, Han and Xiao \cite{YP} is : Are there four limit cycles in all classes 26$-$31?
	
	In this paper, by combining the algorithm of constructing limit cycles by Hofbauer and So to the real root isolation algorithm proposed in \cite{shigen}, four limit cycles are constructed in classes 28,30 and 31 in Theorem 3.1 and Theorem 3.2 of present paper by automatic search. By using the program \texttt{3DLVzd} written by the authors, 6733 examples are searched with randomly chosen interaction matrices. And among these examples, 350, 3, 21, 1 and 1 limit cycles belong to classes 27, 28, 29, 30 and 31, respectively. This gives an affirmative answer to Yu, Han and Xiao's problem \cite{YP}.
	
	Summarizing the results, we have the following main result.
	
	\noindent\bf{Theorem 1.1.} \normalfont\emph{For each class of the six ones  (Classes 26, 27, 28, 29, 30, and 31)  in Zeeman's classification,  there exists a system with four limit cycles.}


	\section{Automated search algorithm}
	In this section, the algorithmic construction method proposed by Hofbauer and So \cite{Ho} and modified by Lu and Luo \cite{Lu02} and Lian, Lu and Luo \cite{lian} is used to search the limit cycles.
	
	Consider a three-dimensional Lotka-Volterra system
	
	\begin{equation}
		\begin{aligned}
			\dot{x}_i=x_i(\sum_{j=1}^{3}a_{ij}(x_j-1)),
		\end{aligned}
	\end{equation}
	where $a_{ij} <0$ and $\mathbf{1}= (1,1,1)$ is the unique positive equilibrium of system (1). 
	

Suppose that matrix $A = (a_{ij})_{3\times3}$ has a real eigenvalue $\lambda $ and a pair of purely imaginary eigenvalues $\pm \omega i (\omega\ne 0).$ Then there is a transformation matrix $T$ to transform A into a block diagonal matrix form $$TAT^{-1}=\left[\begin{array}{ccc}
	c_{11} & c_{12} & 0 
	\\
	c_{21} & c_{22} & 0 
	\\
	0 & 0 & \lambda  
\end{array}\right].
$$
Here, the submatrix has a pair of purely imaginary eigenvalues $\pm\omega i(\omega\ne0)$, that is, the submatrix satisfies $c_{11}+c_{22}= 0 $ and $c_{11}c_{22}-c_{12}c_{21}>0$.
Besides, to guarantee the positive definiteness of the Liapunov function, we need that $c_{21}<0$.

From the center manifold theorem \cite{Carr}, we can suppose that, under the transformation $y=T(x-1)$, the transformed system with linear part $Cy$ has an approximation to the center manifold taking the form

$y_3 = h(y_1,y_2) = h_2(y_1, y_2) + h_3(y_1, y_2) + h_4(y_1, y_2) + h_5(y_1, y_2) + h_6(y_1, y_2) + h.o.t,$\\
where $ y=(y_1,y_2,y_3)^T, h_i(y_1,y_2)=\sum_{j=0}\limits^{k}p_{kj}y_1^{j}y_2^{k-j},(i=2,\dots,6)$ and 
$h.o.t$ denotes the terms with orders greater than or equal to seven.

We have the following steps to calculate the focal value to get the relevant conclusion:

\noindent \begin{tabular}{lp{14.5cm}}
	Step 1:  & Using the subprograms \texttt{randA} and \texttt{DIA (A = randA(), T = DIA(A))}, the coefficient matrix $A$ is  randomly selected, and the block diagonal matrix $TAT^{-1}$ is obtained. \\
	Step 2: &From the subprogram \texttt{Vh}, $c_{ij}$ is solved, so as the approximate center manifold $y_3 = h(y_1, y_2)$ is gotten. \\
	
	Step 3: &Substituting $y_3 = h(y_1, y_2)$ into $\dot{y}_1,\dot{y}_2$ to get a two-dimensional system with center-focus type. \\
	
	Step 4: &Calculating the focal value by the subprogram \texttt{JJLLineSolve}. \\
	
	Step 5: &Applying real root isolation algorithm (subprogram \texttt{mrealroot} \cite{shigen} ) to check the independence of the focal values and get multiple limit cycles.\\
	
	Step 6: &Using subprogram \texttt{CL} to check which class in the Zeeman classification the constructed system belongs to.\\
\end{tabular}

The main algorithm of the program \texttt{3DLVzd} is as follows:
\begin{algorithm}[!ht]
	\KwIn{number}
	\KwOut{coefficient matrix A, Zclass}
	\Begin{
		\For {$zz$ $to$ $numer$}	
		{
			{restart:with(linalg):with(LinearAlgebra):\\}
			{A:=randA(A):P :=DIA(A, $\mu$, [y[1], y[2], y[3]]):\\}
			{g1 := map(factor, P[1]):g2 := map(factor, P[2]):g3 := map(factor, P[3]):\\}
			{hh := factor(Vh2(g1, g2, g3, [y[1], y[2], y[3]], h, 6)):\\}
			{f1 := subs(y[3] = hh, g1):f2 := subs(y[3] = hh, g2):\\}
			{La := factor(JDLLinearSolve(f1, f2, [y[1], y[2]], 3)) :\\}
			{la(i) := expand(numer(La[i])):(i=1,2,3):\\}
			{lad(i) := expand(denom(La[i])):(i=1,2,3):\\}
			{sys :=[la1, la2]:X := [$\lambda$,n]:\\}
			{S := charsets[mcharset](sys, X, basset):\\}
			{S := S[1]:charsets[iniset](S, X):\\}
			{S1 := subs($\lambda$ = -l, S[1][1]):\\}
			{SS := factor(subs(n = -q, S[1][2])):\\}
			{s1 := realroot(S1,$1/10^{10}$):\\}
			{J:=mrealroot(S,{$\lambda$,n}):Check:=CC(A,o,o1):\\}
			{B := subs({n = o,$\lambda$ = o1,$\mu$ = o2}, A):P1 := CL(B):\\}
			{\If{P1[1]=Zclass}
			{\If{$LV_3\in$ number and Check=true}
			 {save A,P1[1],"zd.txt":}}}
		}}
	\caption{\texttt{Automated search algorithm}.\label{alg:1}}
\end{algorithm}
\newpage


Compared with \cite{lian}, it can be seen that the number of terms and degree of the multivariate polynomials increase very rapidly as the order of the focal values increase.

In another aspect,  \cite{lian} used an algorithm for isolating real roots of univariate polynomial systems dealing with focal values, while this paper deals with multivariate polynomial systems based on an algorithm for isolating real roots of multivariate polynomial systems \cite{2005}. In our program, a system of multivariate polynomials is triangularized by using the Epsilon package given by Wang \cite{WWD}. 
\newpage

\section{Three-dimensional Lotka-Volterra competitive systems}

\subsection{Four limit cycles in class 28}
The dynamical behavior of solutions of the class 28 in Zeeman’s classification restricted on the carrying simplex is as follows:
\begin{center}
	\includegraphics[scale=0.1]{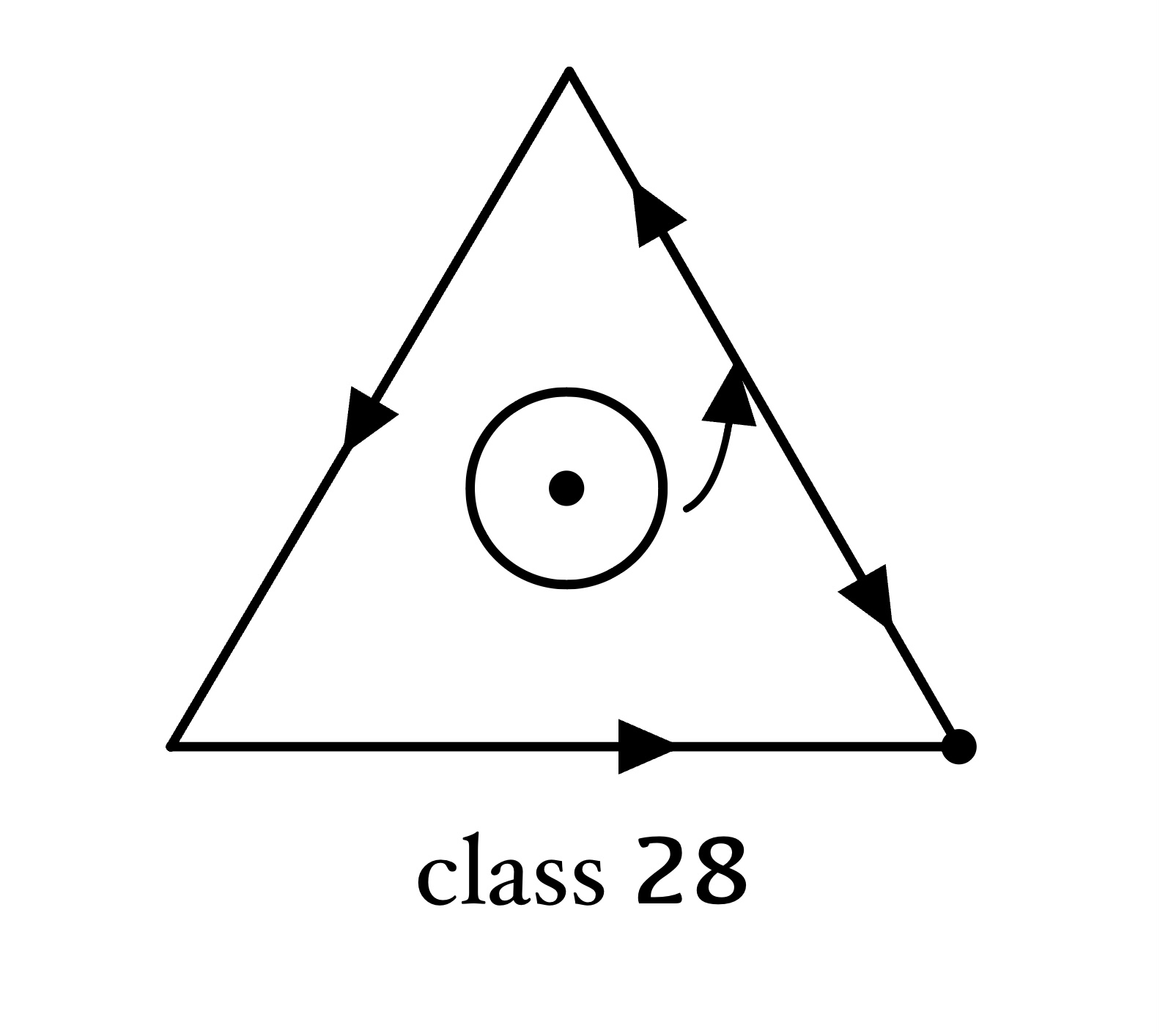}
\end{center}

The found system with an interaction matrix is as follows:

$$A=\left[\begin{array}{ccc}
	-\frac{17}{24} & -2 & -\lambda  
	\\
	-\frac{33}{23} & -10 & -\mu  
	\\
	-n  & -21 n  & -\frac{99}{37} 
\end{array}\right],(n,\mu,\lambda>0).$$

To satisfy the necessary eigenvalue condition \cite{Ho} $\texttt{det}(A) = (A_{11} + A_{22} + A_{33}) \cdot \texttt{trace}(A)$, we need $\mu=-\frac{607835112 \lambda  n -7773334823}{4864016448 n}$. 

By the transformation $y = T (x-1) $,with $T$ as follows,\begin{spacing}{2}$$\left[\begin{array}{ccc}
	-\frac{33}{23} & \frac{3005}{888} & \frac{607835112 \lambda  n -7773334823}{4864016448 n} 
	\\
	-n  & -21 n  & \frac{257}{24} 
	\\
	\frac{175916243}{65729952}+\frac{4791493 \lambda  n}{1825832} & 21 \lambda  n +\frac{257}{12} & \frac{7622102628 \lambda  n +7773334823}{2432008224 n} 
\end{array}\right].
$$
\end{spacing} 
The three-dimensional system is transformed to be a new one whose linear part is in the block diagonal form
\begin{center}\begin{spacing}{2}$\left[\begin{array}{ccc}
			\frac{1014026}{684499} & \frac{416062526328888 \lambda  n +384134040047899}{3329414394639552 n} & 0 
			\\
			-\frac{8896983 n}{684499} & -\frac{1014026}{684499} & 0 
			\\
			0 & 0 & -\frac{11885}{888} 
		\end{array}\right]
		$.\end{spacing}\end{center}
Furthermore, the three-dimensional system can be reduced to a two-dimensional system by the subprogram $Vh$, and
the first three focal values can be obtained by using the subprogram \texttt{JJLLineSolve},
$$\begin{aligned}LV_1=f(\lambda,n)=\frac{f_1(\lambda,n)}{f_2(\lambda,n)},\\
	LV_2=g(\lambda,n)=\frac{g_1(\lambda,n)}{g_2(\lambda,n)},\\
	LV_3=z(\lambda,n)=\frac{z_1(\lambda,n)}{z_2(\lambda,n)},\end{aligned}$$
where
\begin{align*}
	f_1(\lambda,n) = &-684499 (71175864 n \lambda -30452821)\\& (1867427763509790559220459722237440 \lambda^{3} n^{3}\\&-8066558192490463098597559057769472 \lambda^{2} n^{3}\\&-9284198345879360722318571367478464 \lambda^{2} n^{2}\\&+9630923381872490306204845292994048 \lambda \,n^{3}\\&+18212185214244398672238510809517312 \lambda \,n^{2}\\&+5385169712442285368618043473601672 n \lambda \\&+37509227186769280161709353461815488 n^{2}\\&+1036086857152915319628573370644784 n \\&-604462354449619944145534311192809),\\\\
	
	f_2(\lambda,n)=&7401027927168 n^{3} (1169277093792 \lambda  n +31737803345137)\\& (146159636724 \lambda  n +16056506276339)^{2},\\\\
	
	g_2(\lambda,n)=&2717701945572392898631507132761600 (1315436730516 \lambda  n +15556227332951)\\& (1169277093792 \lambda  n +31737803345137)^{3} (146159636724 \lambda  n +16056506276339)^{6} \\&(416062526328888 \lambda  n +384134040047899)^{2} n^{5}

	\\\\
	
	z_2(\lambda,n)=&74656592898569603821499140755938066449603169746217765376000 \\&(4677108375168 \lambda  n +30236966514973) (146159636724 \lambda  n +16056506276339)^{10}\\& (1169277093792 \lambda  n +31737803345137)^{5} (1315436730516 \lambda  n +15556227332951)^{2} \\&(71175864 \lambda  n -30452821) (416062526328888 \lambda  n +384134040047899)^{5} n^{7}
\end{align*}

Here, the lengthy expressions of $g_1(\lambda,n)$ (polynomial of $63$ terms with degree $26$) and $z_1(\lambda,n)$ (polynomial of $169$ terms with degree $50$) are listed in the appendix.

Running the program \texttt{mrealroot} \cite{shigen},

\begin{center}
	>-\texttt{mrealroot}$([f_1(\lambda,n),g_1(\lambda,n)],[\lambda,n],\displaystyle\frac {1}{10^{20}},[z_1(\lambda,n),z_2(\lambda,n),$\texttt{det}(A)]),
\end{center}
we get

\begin{small}
		$$[[\lambda_1,n_1][+,-,-]],$$
		where
		$$\lambda_1\in({\frac{48083713211257141381227}{9444732965739290427392}}, {\frac
			{48083713211499877152963}{9444732965739290427392}}),$$ $$n_1\in({\frac{
				18617876387518095278417715070016705816645420386546066217507077513538675027145
			}{
				57896044618658097711785492504343953926634992332820282019728792003956564819968
		}},$$
	$$ {\frac{
				9308938193759047639208857535008352908322710193273033108753538756769337513573
			}{
				28948022309329048855892746252171976963317496166410141009864396001978282409984
		}}).$$
		
	\end{small}
	
	This shows that it is a competitive system such that the real root of $LV_1 =LV_2=0$ in the interval form is $[\lambda_1,n_1]$ which makes $LV_3<0.$
	
	Thus, we can perturb $LV_1$ and $LV_2$ by using $n$ and $\lambda$ to obtain two small-amplitude limit cycles. Further, we change $\mu$ to perturb $LV_0$ such that $LV_0\cdot LV_1 < 0$ and $|LV_0| \ll |LV_1|$ to get one more small-amplitude limit cycle. These ensure the existence of three small-amplitude limit cycles, and the outer one being stable.
	
	Finally, we check if we are luck to be able to apply Poincaré-Bendixson theorem to get a fourth limit cycle. Using Zeeman’s notation \cite{Zee}, we have we have $R_{ij} =\textrm{sgn}(\alpha_{ij})$  and $Q_{kk}=\textrm{sgn}(\beta_{kk})$, with $\alpha_{ij}=\frac{b_ia_{ji}}{a_{ii}}-b_j=(AR_i)_j-b_j$ and $\beta_{kk}=(AQ_k)_k-b_k$, which are the algebraic invariants of $A$. Here, $R_i$ is the equilibrium on the $x_i$-axis, and $Q_k$ is the positive equilibrium on the plane of $x_k=0$.
	
	Since$$R_{12}=Q_{33}=R_{21}=-R_{23}=R_{32}=-R_{31}=R_{13}=1,$$it implies that the constructed example belongs to class 28 in Zeeman's classification. Since the third focal value is negative, the outer limit cycle is stable. For class 28, the boundary of the simplex is a attractor which can reveal the existence of a fourth limit cycle.  Then by Poincaré-Bendixson theorem we have the following result in this subsection.

\noindent\bf{Theorem 3.1.} \normalfont\emph{There exist at least four limit cycles for class 28 in Zeeman's classification.}

\subsection{Four limit cycles in classes 30 and 31}
The dynamical behavior of solutions of the classes 30 and 31 in Zeeman’s classification restricted on the carrying simplices are as follows:
\begin{center}
	\includegraphics[scale=0.5]{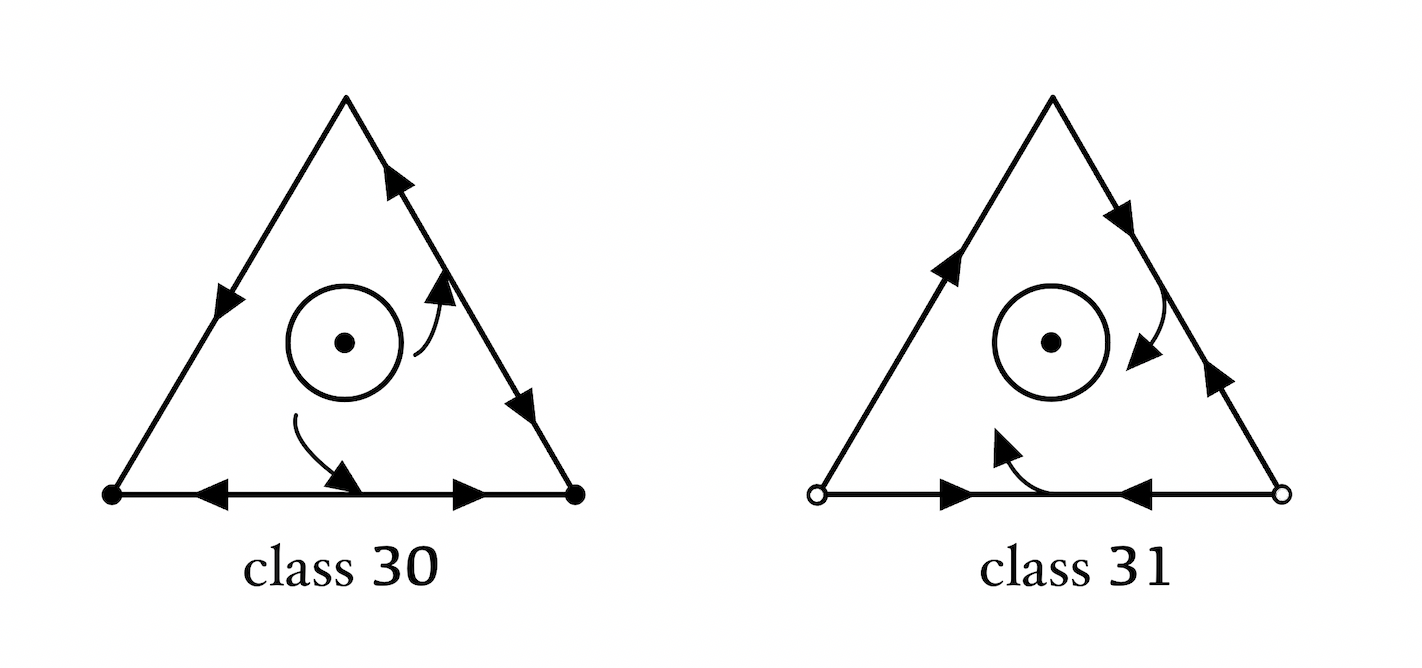}
\end{center}

Obviously, for class 30 when there exist three small-amplitude limit cycles with the outer one being stable and the boundary of the simplex is a attractor,a fourth limit cycle can be obtained based on the Poincaré-Bendxison Theorem. And for class 31 when there exist three small-amplitude limit cycles with the outer one being unstable and the boundary of the simplex is a repellor, a fourth limit cycle can also be obtained.
Among the 6733 examples chosen randomly, there are only one example belonging to class 30 and the other to class 31, respectively , with at least four limit cycles.

 The found system belonging to class 30 with an interaction matrix is as follows:

$$A=\left[\begin{array}{ccc}
	-\frac{31}{26} & -\frac{50}{3} & -\lambda  
	\\
	-\frac{7}{80} & -\frac{10}{9} & -\mu  
	\\
	-n  & -74 n  & -\frac{20}{33} 
\end{array}\right]
,(n,\mu,\lambda>0).$$

To satisfy the necessary eigenvalue condition \cite{Ho} $\texttt{det}(A) = (A_{11} + A_{22} + A_{33}) \cdot \texttt{trace}(A)$, we need $\mu=-\frac{298990926 \lambda  n -135667525}{5194519200 n}$. 

By the transformation $y = T (x-1) $,with $T$ as follows,\begin{spacing}{2}$$\left[\begin{array}{ccc}
		-\frac{7}{80} & \frac{1543}{858} & \frac{298990926 \lambda  n -135667525}{5194519200 n} 
		\\
		-n  & -74 n  & \frac{539}{234} 
		\\
		\frac{76522369}{34630128}+\frac{15758781 \lambda  n}{3699800} & 74 \lambda  n +\frac{13475}{351} & \frac{2876594526 \lambda  n +1492342775}{3428382672 n} 
	\end{array}\right]
	$$
\end{spacing} 
The three-dimensional system is transformed to be a new one whose linear part is in the block diagonal form
\begin{center}\begin{spacing}{2}$\left[\begin{array}{ccc}
			\frac{1014026}{684499} & \frac{416062526328888 \lambda  n +384134040047899}{3329414394639552 n} & 0 
			\\
			-\frac{8896983 n}{684499} & -\frac{1014026}{684499} & 0 
			\\
			0 & 0 & -\frac{11885}{888} 
		\end{array}\right]
		$.\end{spacing}\end{center}
Furthermore, the three-dimensional system can be reduced to a two-dimensional system by the subprogram $Vh$, and
the first three focal values can be obtained by using the subprogram \texttt{JJLLineSolve},
$$\begin{aligned}LV_1=f(\lambda,n)=\frac{f_1(\lambda,n)}{f_2(\lambda,n)},\\
	LV_2=g(\lambda,n)=\frac{g_1(\lambda,n)}{g_2(\lambda,n)},\\
	LV_3=z(\lambda,n)=\frac{z_1(\lambda,n)}{z_2(\lambda,n)},\end{aligned}$$
where
\begin{align*}
	f_1(\lambda,n) = &-141971 \left(2821801554 \lambda  n -580250225\right) (231257128828134319683769384515213372 \lambda^{3} n^{3}\\&-225972234850949262155629372216181760 \lambda^{2} n^{3}\\&+4060846583696338723048867472264040 \lambda^{2} n^{2}\\&+106909009589966694922567133026506000 \lambda  \,n^{3}\\&+94050940211165031622450092875776500 \lambda  \,n^{2}\\&-41106218218019467213434054161387625 \lambda  n \\&-29492401202744021278057001956275000 n^{2}\\&-17295970200081170598648562320656250 n \\&+6047832877840035486968330008071875),\\\\
	
	f_2(\lambda,n)=&1181707100040896 n^{3} (3072941892306 \lambda  n +7349020223275)^{2}\\& (1536470946153 \lambda  n +681667842275)
	,\\\\
	
	g_2(\lambda,n)=&932776164021058410873025903342326686720 \\& (27656477030754 \lambda  n +2293880263075) (1536470946153 \lambda  n +681667842275)^{3} \\& (3072941892306 \lambda  n +7349020223275)^{6}(466928448306606 \lambda  n -72042498603775)^{2} n^{5}

	\\\\
	
	z_2(\lambda,n)=&565016330393399308196222147719410067561468806824224095993338871930880 \\&(12291767569224 \lambda  n -532341800525) (3072941892306 \lambda  n +7349020223275)^{10}\\& (1536470946153 \lambda  n +681667842275)^{5} (27656477030754 \lambda  n +2293880263075)^{2}\\& (2821801554 \lambda  n -580250225) (466928448306606 \lambda  n -72042498603775)^{5} n^{7}
	
\end{align*}

Here, the lengthy expressions of $g_1(\lambda,n)$ (polynomial of $63$ terms with degree $26$) and $z_1(\lambda,n)$ (polynomial of $169$ terms with degree $50$) are listed in the appendix.

Running the program \texttt{mrealroot} \cite{shigen},

\begin{center}
	>-\texttt{mrealroot}$([f_1(\lambda,n),g_1(\lambda,n)],[\lambda,n],\displaystyle\frac {1}{10^{20}},[z_1(\lambda,n),z_2(\lambda,n),$\texttt{det}(A)]),
\end{center}
we get

\begin{small}
		$$[[\lambda_2,n_2][-,+,-]],$$
		where
		$$\lambda_2\in({\frac{684713719509459395319}{590295810358705651712}}, {\frac{
				2738854878037837581277}{2361183241434822606848}}),$$ $$n_2\in({\frac{171259638388598142557329817600}{
				966063118878456319705527262863}},{\frac{209290355979432427921}{1180591620717411303424}}).$$
		
	\end{small}
	
	This shows that it is a competitive system such that the real root of $LV_1 =LV_2=0$ in the interval form is $[\lambda_2,n_2]$ which makes $LV_3<0.$
	
	Thus, we can perturb $LV_1$ and $LV_2$ by using $n$ and $\lambda$ to obtain two small-amplitude limit cycles. Further, we change $\mu$ to perturb $LV_0$ such that $LV_0\cdot LV_1 < 0$ and $|LV_0| \ll |LV_1|$ to get one more small-amplitude limit cycle. These ensure the existence of three small-amplitude limit cycles, and the outer one being stable.
	
	Similarly, we need to check whether we are able to apply Poincaré-Bendixson theorem to obtain a fourth limit cycle. Using Zeeman’s notation \cite{Zee}, we have $$R_{12}=Q_{33}=R_{21}=R_{23}=R_{32}=Q_{11}=-R_{31}=R_{13}=1,$$it implies that the constructed example belongs to class 30 in Zeeman's classification. Since the third focal value is negative, the outer limit cycle is stable. For class 30, the boundary of the simplex is a attractor which can reveal the existence of a fourth limit cycle. 

For an example belonging to class 31, the found system with an interaction matrix is as follows:
\begin{spacing}{2}
	$$A=\left[\begin{array}{ccc}
		-\frac{17}{2} & -\frac{25}{3} & -\frac{5 \lambda}{6} 
		\\
		-\frac{1}{2} & -\frac{19}{2} & -\frac{5 \mu}{6} 
		\\
		-\frac{5 n}{6} & -15 n  & -\frac{31}{6} 
	\end{array}\right], (n , \mu , \lambda>0).$$
\end{spacing}
\begin{spacing}{2}To satisfy the necessary eigenvalue condition \cite{Ho} $\texttt{det}(A) = (A_{11} + A_{22} + A_{33}) \cdot \texttt{trace}(A)$, we need $\mu=-\frac{4 \left(425 n \lambda -95391\right)}{20425 n}$. 
	
	By the transformation $y = T (x-1) $,with $T$ as follows,\end{spacing} 
\begin{center}\begin{spacing}{2}$\left[\begin{array}{ccc}
			-\frac{1}{2} & \frac{41}{3} & \frac{850 n \lambda -190782}{12255 n} 
			\\
			-\frac{5 n}{6} & -15 n  & 18 
			\\
			\frac{10200}{817}+\frac{850 n \lambda}{817} & \frac{25 n \lambda}{2}+150 & \frac{17665 n \lambda +211980}{1634 n} 
		\end{array}\right]
		,$\end{spacing}\end{center}
	
The three-dimensional system is transformed to be a new one whose linear part is in the block diagonal form
\begin{center}\begin{spacing}{2}$\left[\begin{array}{ccc}
			-\frac{475}{136} & \frac{115600 \lambda  n -17286969}{1666680 n} & 0 
			\\
			-\frac{2035 n}{408} & \frac{475}{136} & 0 
			\\
			0 & 0 & -\frac{139}{6} 
		\end{array}\right]
		$.\end{spacing}\end{center}
Furthermore, the three-dimensional system can be reduced to a two-dimensional system by the subprogram $Vh$, and
the first three focal values can be obtained by using the subprogram \texttt{JJLLineSolve},
$$\begin{aligned}LV_1=f(\lambda,n)=\frac{f_1(\lambda,n)}{f_2(\lambda,n)},\\
	LV_2=g(\lambda,n)=\frac{g_1(\lambda,n)}{g_2(\lambda,n)},\\
	LV_3=z(\lambda,n)=\frac{z_1(\lambda,n)}{z_2(\lambda,n)},\end{aligned}$$
where
\begin{align*}
	f_1(\lambda,n) = &-272 (10175 n \lambda -1880367) (1615879075312500 \lambda^{3} n^{3}\\&+5102402771778125 \lambda^{2} n^{3}-318360822448914375 \lambda^{2} n^{2}\\&+23010457575345625 \lambda \,n^{3}+1091326301880441125 \lambda \,n^{2}\\&+4265788315448606700 n \lambda +4943072207163509100 n^{2}\\&-53668496608032586500 n +29081719740318720336),\\\\
	
	f_2(\lambda,n)=&1552959375 n^{3} (40700 n \lambda +8263789) (2035 n \lambda +2780978)^{2},\\\\
	
	g_2(\lambda,n)=&378989256106975830078125 n^{5} (2475 n \lambda -30758) (115600 n \lambda -17286969)^{2}\\& (40700 n \lambda +8263789)^{3} (2035 n \lambda +2780978)^{6},\\\\
	
	z_2(\lambda,n)=&158102013753180527268510036960601806640625 n^{7} (32560 n \lambda -2860123) \\&(10175 n \lambda -1880367) (2475 n \lambda -30758)^{2} (40700 n \lambda +8263789)^{5}\\& (115600 n \lambda -17286969)^{5} (2035 n \lambda +2780978)^{10}.

\end{align*}
Here, the lengthy expressions of $g_1(\lambda,n)$ (polynomial of $63$ terms with degree $26$) and $z_1(\lambda,n)$ (polynomial of $169$ terms with degree $50$) are listed in the appendix.

Running the program \texttt{mrealroot} \cite{shigen},

\begin{center}
	>-\texttt{mrealroot}$([f_1(\lambda,n),g_1(\lambda,n)],[\lambda,n],\displaystyle\frac {1}{10^{10}},[z_1(\lambda,n),z_2(\lambda,n),$\texttt{det}(A)]),
\end{center}
we get

\begin{small}
		$$[[\lambda_3,n_3][+,+,-]],$$
		where
		$$\lambda_3\in(\frac{161866876864428888044489}{2361183241434822606848}, \frac{80933438432214444022245}{1180591620717411303424}),$$ $$n_3\in({\frac{4439891048147073080770953216}{
				1646995472095563935852675575}}, {\frac{3182582015038192507591}{
				1180591620717411303424}}).$$
		
	\end{small}
	
	This shows that it is a competitive system such that the real root of $LV_1 =LV_2=0$ in the interval form is $[\lambda_3,n_3]$ which makes $LV_3>0.$
	
	
	
	Thus, we can perturb $LV_1$ and $LV_2$ by using $n$ and $\lambda$ to obtain two small-amplitude limit cycles. Further, we change $\mu$ to perturb $LV_0$ such that $LV_0\cdot LV_1 < 0$ and $|LV_0| \ll |LV_1|$ to get one more small-amplitude limit cycle. These ensure the existence of three small-amplitude limit cycles, and the outer one being unstable.
	
	Similarly, we check if we are able to apply Poincaré-Bendixson theorem to obtain a fourth limit cycle. Using Zeeman’s notation \cite{Zee}, we have $$R_{12}=Q_{33}=R_{21}=R_{23}=R_{32}=Q_{11}=-R_{31}=R_{13}=-1,$$it implies that the constructed example belongs to class 31 in Zeeman's classification. Since the third focal value is positive, the outer limit cycle is unstable. For class 31, the boundary of the simplex is a repellor which can reveal the existence of a fourth limit cycle.  
	Then by Poincaré-Bendixson theorem we have the following result in this subsection.

	\noindent\bf{Theorem 3.2.} \normalfont\emph{There exist at least four limit cycles for class 30 and class 31 in Zeeman's classification.}
	
	\section{Concluding remarks}
	In this paper, four limit cycles are constructed in classes 28, 30 and 31 of Zeeman's classification. By incorporating Hirsch's dimension reduction method, the center manifold construction, the focal values computation   combined with the real root isolation algorithm  with multivariable polynomials \cite{2001,2005} and automatic search, the independence of the focal values is shown  and the same stability of the  outermost small amplitude limit cycle and  the boundary of the simplex are checked automatically, thus three small amplitude limit cycles  and a large scale  limit cycle are obtained. Therefore, in each of six classes of Zeeman's classification, the existence of four limit cycles is ensured.
	
	Obviously, how to construct five limit cycles for systems in Zeeman's classification of six types (classes 26$-$31) is an interesting and challenging problem. After obtaining the focal values up to order 4, the remaining problem is to triangulate these large polynomials so as to check their independence and to ensure the same stability of the outermost  small amplitude limit cycle and the boundary of the simplex. To obtain  directly focal values up to 5th order and prove their independence may be more complicated.


\section*{Appendix}
The above $g_1(\lambda,n),z_1(\lambda,n)$ are listed in the https://sourl.cn/VZN2nR \href{https://sourl.cn/VZN2nR}.
	
\end{document}